\documentclass[12pt]{article}

\usepackage[utf8]{inputenc}
\usepackage[russian,english]{babel}
\usepackage{latexsym}
\usepackage{amsfonts,amsmath,amssymb,amsthm,mathtools} 
\usepackage{amscd}
\usepackage{graphicx}
\usepackage{multicol}		% набор в несколько колонок
\usepackage[dvipsnames]{color}	% цветной текс
\usepackage{mathrsfs}		% пакеты специальных шрифтов
\usepackage{longtable}
\usepackage{lscape}
\theoremstyle{remark}
\newtheorem*{remark}{Remark}
\theoremstyle{plain}
\newtheorem*{statement}{Statement}

\DeclareMathOperator{\im}{im}
\DeclareMathOperator{\rank}{rank}

\mathtoolsset{showonlyrefs=true}	% показывать только те номера формул на
\author{Nurlan M. Sadykov}
\title{Two-color solutions of set-theoretic tetrahedron equation and their cohomologies}
\begin{document}

\maketitle

\begin{abstract}
All solutions of the set-theoretic constant tetrahedron equation with two colors
are found, and some of their properties are analyzed. The list includes 406
solutions --- we call them $R$-operators,--- most of which are degenerate
(non-bijective). Then, we calculate the 3-cohomologies for our $R$-operators,
and discuss the applicability of our results to 3-dimensional statistical
physics.
\end{abstract}

\section{Introduction}

\subsection{Set-theoretic tetrahedron equation}
Consider some set $X$ and mapping $R$
\begin{equation}
	X\times X\times X \xrightarrow{R} X \times X\times X
\end{equation}
We define an operator $R_{ijk}\colon\;X^{\times 6} \to X^{\times 6}$ as
mapping~$R$ acting on the $i$th, $j$th and $k$th factors in~$X^{\times 6}$, for
instance, if $R(x_1,x_2,x_3)=(x_1',x_2',x_3')$, then
\[
R_{123}(x_1,x_2,x_3,x_4,x_5,x_6)=(x_1',x_2',x_3',x_4,x_5,x_6).
\]
We will consider $R$-operators that are solutions of the following
\emph{set-theoretic tetrahedron equation} (STTE):
\begin{equation}
	\label{eq:tq}
	R_{123}R_{145}R_{246}R_{356} = R_{356}R_{246}R_{145}R_{123}.
\end{equation}

We can construct a solution of the quantum version of tetrahedron
equation from a solution of STTE. Let $V$ be a linear space over some field~$\mathbb K$
(or maybe even a module over a commutative ring) with a basis $\{e_x\}$ in~$V$
numbered with index $x\in V$. We
define the quantum \emph{linear} operator $\mathcal{R}\colon\; V^{\otimes 3} \to V^{\otimes 3}$ as
follows:
\begin{equation}
	\label{eq:R-q-prostoe}
	\mathcal{R}(e_{x_1} \otimes e_{x_2} \otimes e_{x_3}) = 
	e_{r_1(x_1, x_2, x_3)} \otimes e_{r_2(x_1,x_2,x_3)} \otimes
	e_{r_3(x_1,x_2,x_3)},
\end{equation}
where $r_1$, $r_2$ and~$r_3$ are, of course, the three components of mapping~$R$.

Such solutions have been studied by Hietarinta~\cite{Hietarinta}.

\subsection{Tetrahedral chain complex and related definitions}\label{ss:complex}
The \emph{tetrahedral chain complex} was introduced in paper~\cite{KST}. Here are the definitions that
we will need, and some related comments.

In an $n$-dimensional cube~$I^n$, we call \emph{incoming} each of its $(n-1)$-faces
in the following sequence:
\[\{0\}\times I\times I \times\dots,\quad I\times\{1\}\times I\times \dots,\quad I\times I\times\{0\}\times\dots,\quad \dots ,\]
and the opposite faces
\[\{1\}\times I\times I \times\dots,\quad I\times\{0\}\times I\times \dots,\quad I\times I\times\{1\}\times\dots,\quad \dots ,\]
are called \emph{outgoing}.

A \emph{2-coloring} of~$I^n$ is a map~$\vartheta$ from the set of its 2-faces into our set~$X$.
Set~$X$ is called in this context the \emph{set of colors}.

A coloring of a \emph{three-dimensional} cube is
called \emph{permitted} if the colors of its three outgoing 2-faces are obtained from
the colors of its three incoming 2-faces via a map $R$, with $R$ being a (fixed) solution
of STTE on $X$. A coloring $\vartheta$ of $n$-dimensional cube is permitted if its
restriction on any 3-subcube is permitted.

The space~$C^n(X)$ of \emph{$n$-chains}
is the free $\mathbb{Z}$-module spanned by the set of all permitted
colorings of the cube $I^n$.

\begin{remark}
In this paper, we consider thus only chains with integer coefficients, but our \emph{cochains}
will sometimes take values in a multiplicative group of a field, see below.
\end{remark}

The \emph{boundary operator}~$\partial$ is defined on permitted colorings (which form, as we have
just said, a $\mathbb{Z}$-basis for all chains) in the following way:
\begin{equation}\label{eq:border}
	\partial( \vartheta ) = \sum_{k=1}^n(d_k^{\rm in}\vartheta -
	d_k^{\rm out}\vartheta),
\end{equation}
where $d^{\rm in}_k$ is the restriction of coloring $\vartheta$ on the $k$th
$(n-1)$-dimensional incoming cube, and $d^{\rm out}_k$ is the restriction on the $k$th outgoing cube.

The cochains and the corresponding coboundary operator are defined in the standard dual way.
Although we will be dealing mainly with $\mathbb Z$-cochains, we would like to demonstrate
how the above definition~\eqref{eq:R-q-prostoe} of quantum~$\mathcal R$ can be modified, using a
3-cocycle taking values in the \emph{multiplicative} group~$\mathbb K^*$ of field~$\mathbb K$.

Let $c(x_1,x_2,x_3)$ be a 3-cocycle (to be exact, $c(x_1,x_2,x_3)$ is of course its restriction
onto the basis elements)
of tetrahedral chain complex taking values in~$\mathbb K^*$.  Then, it can be checked that the quantum
operator defined as follows
\begin{equation}
	\label{eq:Roperator}
	\mathcal{R}(e_{x_1} \otimes e_{x_2} \otimes e_{x_3}) = c(x_1,x_2,x_3)
	e_{r_1(x_1, x_2, x_3)} \otimes e_{r_2(x_1,x_2,x_3)} \otimes
	e_{r_3(x_1,x_2,x_3)}
\end{equation}
still satisfies the quantum tetrahedron equation.

\subsection{The results of this paper}
In this paper we find all solutions of STTE for the case of two-element set $X$,
$\# X = 2$. Moreover, we calculate integer 3-cohomologies for all of them.

We write our $R$-operator in components~$r_k$ as in \eqref{eq:R-q-prostoe} and~\eqref{eq:Roperator}: 
\[
R(x,y,z) = (r_1(x,y,z), r_2(x,y,z), r_3(x,y,z)).
\]
In the case of two-element set~$X$, we
can identify~$X$ with~$\mathbb{Z}_2$, so $r_k\colon\;\mathbb{Z}_2^{\times 3} \to
\mathbb{Z}_2$. Moreover, each~$r_k$ can be written as a polynomial in $x$, $y$
and~$z$ of degree at most one in each variable.

We have used a brute force method of searching for solutions of STTE through all $2^{24}$ triples $(r_1,r_2,r_3)$.

Below,
\begin{itemize}\itemsep 0pt
 \item in Section~\ref{s:explanation}, we explain the general structure of (the presentation of) our results,
 \item in the main Section~\ref{s:main}, the catalogue of all 406 obtained $R$-operators is presented, and
 \item in Section~\ref{s:discussion}, we discuss our results, explaining, in particular, how \emph{degenerate} $R$-operators may be applied to statistical physics.
\end{itemize}

\section{What exactly is calculated, and how the results are presented}\label{s:explanation}

The (hopefully easy to use) catalogue of our results is presented below in Section~\ref{s:main}.
Here we explain how our results are presented and what exactly they include.

\subsection{Two simple symmetries}\label{ss:sym}
Consider two transformations that take a triple and return a triple:
\begin{equation}
	\sigma_1: (x,y,z) \mapsto (z,y,x),
\end{equation}
and
\begin{equation}
	\sigma_2: (x,y,z) \mapsto (1-x, 1-y, 1-z).
\end{equation}
We can see that if $R$ is a solution of STTE, then map $R'=\sigma_l R \sigma_l$, where
$l=1,2$, is a solution, too.

\subsection{Cohomologies}

We have calculated 3-cohomologies over the ring~$\mathbb Z$ for the complex introduced 
in Subsection~\ref{ss:complex}. These cohomologies are the factor $\ker\delta^3 /
\im\delta^2$, where $\delta$ is the coboundary operator.

As the basis in the space of 3-cochains consists of triples $(a,b,c) \in
	\mathbb{Z}_2^3$, we represent any vector in this space as a column of length eight,
using the \emph{lexicographic} ordering of basis vectors: $(0,0,0)$, $(0,0,1)$, \dots, 
$(1,1,1)$.

Similarly, we identify $k$-cochains with column vectors for any~$k$. And, accordingly,
we identify the coboundary operators with matrices acting on these vectors from the left.

As for $k$-\emph{chains}, we represent them below as \emph{row} vectors,
and the \emph{boundary} operators acting on them will be identified with the \emph{same} matrices
as for coboundary operators, but acting on rows from the \emph{right}.

An important comment is that the cochain
\begin{equation}
        \label{eq:v_1}
		v_1=(1,1,1,1,1,1,1,1)^{\rm T}
\end{equation}
is \emph{always} 
in~$\ker\delta^3$, making thus a trivial part of cohomologies (and, actually, suggesting 
a change in the definition of our chain complex --- but here we follow the [first three
arXiv versions of] paper~\cite{KST}).

\begin{statement}
For a finite set~$X$ of colors,
\[\dim \im \delta^2 < \# X.\]
\end{statement}

\begin{proof}
Consider a basis vector $v \in X^{\times 3}$ corresponding to a permitted coloring
of the 3-cube. In the cube, there are three incoming and three outgoing 2-faces.
The boundary operator~$\partial^3$ sends~$v$ to a 2-chain that consists just of elements
in~$X$ with integer coefficients. The coefficient at some element $x\in X$ is made 
by adding up values $\pm 1$ with the sign depending on whether $v$ has $x$ as the
color of an incoming or outgoing 2-face.

As the number of incoming faces is the same as of outgoing ones, the sum of all
coefficients in the 2-chain $\partial^3 v$ vanishes. For matrix $\partial^3$,
	which, as we agreed, acts on rows from the right, 
this means that the sum of its entries in any row vanishes, so the columns
of~$\partial^3$
are linearly dependent, and $\rank \partial^2$ is less than the number of columns,
i.e., $\# X$.

	Now, the statement follows from the fact that the coboundary matrix coincides,
        according to our agreement, with the corresponding boundary matrix: \[\partial^3 = \delta^2.\]
\end{proof}

\begin{remark}
The above Statement can be easily generalized onto other dimensions. For instance,
	$\rank\partial^4 < (\# X)^3.$
\end{remark}

In our case $\#X=2$, our calculations show that, indeed, always $\dim \im \delta^2 \le 1$, and typically $\dim \im \delta^2 = 1$.

\subsection{Organization of the catalogue}

The found solutions $R=R_1,\dots,R_{406}$ are sorted in the increasing order of cardinality of the
image of map~$R$. Accordingly, Section~\ref{s:main} is subdivided into eight
subsections.

In the case if the kernel $\ker\delta^3$ does \emph{not} coincide with its 
trivial part, that is, the $\mathbb Z$-linear span of $\im
\delta^2$ and vector $v_1$ given by~\eqref{eq:v_1},we write out the matrices
\emph{transposed} (for convenience) to matrices whose rows $\mathbb{Z}$-span
the mentioned image and kernel.
In addition, all solutions are grouped into subsets of one, two or four
$R$-operators that are taken into each other by symmetries $\sigma_1$ and
$\sigma_2$ defined in Subsection~\ref{ss:sym}. The symmetries are shown in the
form of diagrams. If $R$ is self-symmetric w.r.t\ one of $\sigma_i$, this~$\sigma_i$
is not shown, and if $R$ is self-symmetric w.r.t\ both symmetries,
we write the empty set symbol instead of a diagram.

\section{The results of calculation}\label{s:main}

% версия вывода от 05.04.2015
\begin{enumerate}

\subsection{Solutions of image cardinality 1}
\item \begin{equation}
\begin{CD}
R_{1}\\
 @V\sigma_2VV \\
R_{2}
\end{CD}
\quad\quad\quad
\begin{aligned}
R_{1} &=	( 0,\ 0,\ 0) \\
R_{2} &=	( 1,\ 1,\ 1) 
\end{aligned}
\end{equation}

\subsection{Solutions of image cardinality 2}
\item \begin{equation}
\begin{CD}
R_{3} @>\sigma_1>> R_{4}\\
 @V\sigma_2VV 	 @VV\sigma_2V \\
R_{5} @>>\sigma_1> R_{6}
\end{CD}
\quad\quad\quad
\begin{aligned}
R_{3} &=	( 0,\ 0,\ y) \\
R_{4} &=	( y,\ 0,\ 0) \\
R_{5} &=	( 1,\ 1,\ y) \\
R_{6} &=	( y,\ 1,\ 1) 
\end{aligned}
\end{equation}
\item \begin{equation}
\begin{CD}
R_{7} @>\sigma_1>> R_{8}\\
 @V\sigma_2VV 	 @VV\sigma_2V \\
R_{9} @>>\sigma_1> R_{10}
\end{CD}
\quad\quad\quad
\begin{aligned}
R_{10} &=	( 0,\ 0,\ z) \\
R_{7} &=	( x,\ 0,\ 0) \\
R_{8} &=	( 1,\ 1,\ z) \\
R_{9} &=	( x,\ 1,\ 1) 
\end{aligned}
\end{equation}
\item \begin{equation}
\begin{CD}
R_{11} @>\sigma_1>> R_{12}\\
 @V\sigma_2VV 	 @VV\sigma_2V \\
R_{13} @>>\sigma_1> R_{14}
\end{CD}
\quad\quad\quad
\begin{aligned}
R_{11} &=	( 0,\ 0,\ y+z) \\
R_{12} &=	( x+y,\ 0,\ 0) \\
R_{13} &=	( 1,\ 1,\ y+z+1) \\
R_{14} &=	( x+y+1,\ 1,\ 1) 
\end{aligned}
\end{equation}
\item \begin{equation}
\begin{CD}
R_{15} @>\sigma_1>> R_{16}\\
 @V\sigma_2VV 	 @VV\sigma_2V \\
R_{17} @>>\sigma_1> R_{18}
\end{CD}
\quad\quad\quad
\begin{aligned}
R_{15} &=	( 0,\ 0,\ yz+y+z) \\
R_{16} &=	( xy+x+y,\ 0,\ 0) \\
R_{17} &=	( 1,\ 1,\ yz) \\
R_{18} &=	( xy,\ 1,\ 1) 
\end{aligned}
\end{equation}
\item \begin{equation}
\begin{CD}
R_{19} @>\sigma_1>> R_{20}\\
 @V\sigma_2VV 	 @VV\sigma_2V \\
R_{21} @>>\sigma_1> R_{22}
\end{CD}
\quad\quad\quad
\begin{aligned}
R_{19} &=	( 0,\ x,\ 0) \\
R_{20} &=	( 0,\ z,\ 0) \\
R_{21} &=	( 1,\ x,\ 1) \\
R_{22} &=	( 1,\ z,\ 1) 
\end{aligned}
\end{equation}
\item \begin{equation}
\begin{CD}
R_{23}\\
 @V\sigma_2VV \\
R_{24}
\end{CD}
\quad\quad\quad
\begin{aligned}
R_{23} &=	( 0,\ y,\ 0) \\
R_{24} &=	( 1,\ y,\ 1) 
\end{aligned}
\end{equation}
\item \begin{equation}
\begin{CD}
R_{25} @>\sigma_1>> R_{26}\\
 @V\sigma_2VV 	 @VV\sigma_2V \\
R_{27} @>>\sigma_1> R_{28}
\end{CD}
\quad\quad\quad
\begin{aligned}
R_{25} &=	( 0,\ y,\ y) \\
R_{26} &=	( y,\ y,\ 0) \\
R_{27} &=	( 1,\ y,\ y) \\
R_{28} &=	( y,\ y,\ 1) 
\end{aligned}
\end{equation}
\item \begin{equation}
\begin{CD}
R_{29} @>\sigma_1>> R_{30}\\
 @V\sigma_2VV 	 @VV\sigma_2V \\
R_{31} @>>\sigma_1> R_{32}
\end{CD}
\quad\quad\quad
\begin{aligned}
R_{29} &=	( 0,\ x+y,\ 0) \\
R_{30} &=	( 0,\ y+z,\ 0) \\
R_{31} &=	( 1,\ x+y+1,\ 1) \\
R_{32} &=	( 1,\ y+z+1,\ 1) 
\end{aligned}
\end{equation}
\item \begin{equation}
\begin{CD}
R_{33} @>\sigma_1>> R_{34}\\
 @V\sigma_2VV 	 @VV\sigma_2V \\
R_{35} @>>\sigma_1> R_{36}
\end{CD}
\quad\quad\quad
\begin{aligned}
R_{33} &=	( 0,\ z,\ z) \\
R_{34} &=	( x,\ x,\ 0) \\
R_{35} &=	( 1,\ z,\ z) \\
R_{36} &=	( x,\ x,\ 1) 
\end{aligned}
\end{equation}
\item \begin{equation}
\begin{CD}
R_{37}\\
 @V\sigma_2VV \\
R_{38}
\end{CD}
\quad\quad\quad
\begin{aligned}
R_{37} &=	( 0,\ x+z,\ 0) \\
R_{38} &=	( 1,\ x+z+1,\ 1) 
\end{aligned}
\end{equation}
\begin{multline}
{\bf R_{37}\colon}\quad \im\delta^2=
	\begin{pmatrix}
0 & 0 & -1 & -1 & 0 & -2 & -1 & -3 
	\end{pmatrix}^{\rm T}
\\ \ker\delta^3=
	\begin{pmatrix}
1 & 1 & 0 & 0 & 1 & 0 & 0 & -1 \\
0 & 0 & 1 & 1 & 0 & 0 & 1 & 1 \\
0 & 0 & 0 & 0 & 0 & 1 & 0 & 1 
	\end{pmatrix}^{\rm T}
\end{multline}
\begin{multline}
{\bf R_{38}\colon}\quad \im\delta^2=
	\begin{pmatrix}
3 & 1 & 2 & 0 & 1 & 1 & 0 & 0 
	\end{pmatrix}^{\rm T}
\\ \ker\delta^3=
	\begin{pmatrix}
1 & 0 & 0 & -1 & 0 & 0 & -1 & -1 \\
0 & 1 & 0 & 1 & 1 & 1 & 1 & 1 \\
0 & 0 & 1 & 1 & 0 & 0 & 1 & 1 
	\end{pmatrix}^{\rm T}
\end{multline}
\item \begin{equation}
\begin{CD}
R_{39}\\
 @V\sigma_2VV \\
R_{40}
\end{CD}
\quad\quad\quad
\begin{aligned}
R_{39} &=	( 0,\ x+y+z,\ 0) \\
R_{40} &=	( 1,\ x+y+z,\ 1) 
\end{aligned}
\end{equation}
\begin{multline}
{\bf R_{39}\colon}\quad \im\delta^2=
	\begin{pmatrix}
0 & 0 & 0 & -2 & 0 & -2 & -2 & -2 
	\end{pmatrix}^{\rm T}
\\ \ker\delta^3=
	\begin{pmatrix}
1 & 1 & 1 & 0 & 1 & 0 & 0 & 0 \\
0 & 0 & 0 & 1 & 0 & 1 & 1 & 1 
	\end{pmatrix}^{\rm T}
\end{multline}
\begin{multline}
{\bf R_{40}\colon}\quad \im\delta^2=
	\begin{pmatrix}
2 & 2 & 2 & 0 & 2 & 0 & 0 & 0 
	\end{pmatrix}^{\rm T}
\\ \ker\delta^3=
	\begin{pmatrix}
1 & 1 & 1 & 0 & 1 & 0 & 0 & 0 \\
0 & 0 & 0 & 1 & 0 & 1 & 1 & 1 
	\end{pmatrix}^{\rm T}
\end{multline}
\item \begin{equation}
\begin{CD}
R_{41} @>\sigma_1>> R_{42}\\
 @V\sigma_2VV 	 @VV\sigma_2V \\
R_{43} @>>\sigma_1> R_{44}
\end{CD}
\quad\quad\quad
\begin{aligned}
R_{41} &=	( 0,\ xy+x+y,\ 0) \\
R_{42} &=	( 0,\ yz+y+z,\ 0) \\
R_{43} &=	( 1,\ xy,\ 1) \\
R_{44} &=	( 1,\ yz,\ 1) 
\end{aligned}
\end{equation}
\item \begin{equation}
\begin{CD}
R_{45}\\
 @V\sigma_2VV \\
R_{46}
\end{CD}
\quad\quad\quad
\begin{aligned}
R_{45} &=	( 0,\ xz,\ 0) \\
R_{46} &=	( 1,\ xz+x+z,\ 1) 
\end{aligned}
\end{equation}
\item \begin{equation}
\begin{CD}
R_{47}\\
 @V\sigma_2VV \\
R_{48}
\end{CD}
\quad\quad\quad
\begin{aligned}
R_{47} &=	( 0,\ xz+x+z,\ 0) \\
R_{48} &=	( 1,\ xz,\ 1) 
\end{aligned}
\end{equation}
\begin{multline}
{\bf R_{47}\colon}\quad \im\delta^2=
	\begin{pmatrix}
0 & 0 & -1 & -1 & 0 & -1 & -1 & -2 
	\end{pmatrix}^{\rm T}
\\ \ker\delta^3=
	\begin{pmatrix}
1 & 1 & 0 & 0 & 1 & 0 & 0 & -1 \\
0 & 0 & 1 & 1 & 0 & 0 & 1 & 1 \\
0 & 0 & 0 & 0 & 0 & 1 & 0 & 1 
	\end{pmatrix}^{\rm T}
\end{multline}
\begin{multline}
{\bf R_{48}\colon}\quad \im\delta^2=
	\begin{pmatrix}
2 & 1 & 1 & 0 & 1 & 1 & 0 & 0 
	\end{pmatrix}^{\rm T}
\\ \ker\delta^3=
	\begin{pmatrix}
1 & 0 & 0 & -1 & 0 & 0 & -1 & -1 \\
0 & 1 & 0 & 1 & 1 & 1 & 1 & 1 \\
0 & 0 & 1 & 1 & 0 & 0 & 1 & 1 
	\end{pmatrix}^{\rm T}
\end{multline}
\item \begin{equation}
\begin{CD}
R_{49} @>\sigma_1>> R_{50}\\
 @V\sigma_2VV 	 @VV\sigma_2V \\
R_{51} @>>\sigma_1> R_{52}
\end{CD}
\quad\quad\quad
\begin{aligned}
R_{49} &=	( 0,\ yz,\ yz) \\
R_{50} &=	( xy,\ xy,\ 0) \\
R_{51} &=	( 1,\ yz+y+z,\ yz+y+z) \\
R_{52} &=	( xy+x+y,\ xy+x+y,\ 1) 
\end{aligned}
\end{equation}
\item \begin{equation}
\begin{CD}
R_{53} @>\sigma_1>> R_{54}\\
 @V\sigma_2VV 	 @VV\sigma_2V \\
R_{55} @>>\sigma_1> R_{56}
\end{CD}
\quad\quad\quad
\begin{aligned}
R_{53} &=	( 0,\ yz+y+z,\ yz+y+z) \\
R_{54} &=	( xy+x+y,\ xy+x+y,\ 0) \\
R_{55} &=	( 1,\ yz,\ yz) \\
R_{56} &=	( xy,\ xy,\ 1) 
\end{aligned}
\end{equation}
\item \begin{equation}
\begin{CD}
R_{57}\\
 @V\sigma_2VV \\
R_{58}
\end{CD}
\quad\quad\quad
\begin{aligned}
R_{57} &=	( 0,\ xyz,\ 0) \\
R_{58} &=	( 1,\ xyz+xy+xz+yz+x+y+z,\ 1) 
\end{aligned}
\end{equation}
\item \begin{equation}
\begin{CD}
R_{59}\\
 @V\sigma_2VV \\
R_{60}
\end{CD}
\quad\quad\quad
\begin{aligned}
R_{59} &=	( 0,\ xyz+xy+xz+yz+x+y+z,\ 0) \\
R_{60} &=	( 1,\ xyz,\ 1) 
\end{aligned}
\end{equation}
\item \begin{equation}
\begin{CD}
R_{61}\\
 @V\sigma_2VV \\
R_{62}
\end{CD}
\quad\quad\quad
\begin{aligned}
R_{61} &=	( y,\ 0,\ y) \\
R_{62} &=	( y,\ 1,\ y) 
\end{aligned}
\end{equation}
\item \begin{equation}
\varnothing
\quad\quad\quad
\begin{aligned}
R_{63} &=	( y,\ y,\ y) 
\end{aligned}
\end{equation}
\item \begin{equation}
\varnothing
\quad\quad\quad
\begin{aligned}
R_{64} &=	( y,\ y+1,\ y) 
\end{aligned}
\end{equation}

\subsection{Solutions of image cardinality 3}
\item \begin{equation}
\begin{CD}
R_{65} @>\sigma_1>> R_{66}\\
 @V\sigma_2VV 	 @VV\sigma_2V \\
R_{67} @>>\sigma_1> R_{68}
\end{CD}
\quad\quad\quad
\begin{aligned}
R_{65} &=	( 0,\ y,\ yz) \\
R_{66} &=	( xy,\ y,\ 0) \\
R_{67} &=	( 1,\ y,\ yz+y+z) \\
R_{68} &=	( xy+x+y,\ y,\ 1) 
\end{aligned}
\end{equation}
\item \begin{equation}
\begin{CD}
R_{69} @>\sigma_1>> R_{70}\\
 @V\sigma_2VV 	 @VV\sigma_2V \\
R_{71} @>>\sigma_1> R_{72}
\end{CD}
\quad\quad\quad
\begin{aligned}
R_{69} &=	( 0,\ y,\ yz+z) \\
R_{70} &=	( xy+x,\ y,\ 0) \\
R_{71} &=	( 1,\ y,\ yz+y+1) \\
R_{72} &=	( xy+y+1,\ y,\ 1) 
\end{aligned}
\end{equation}
\item \begin{equation}
\begin{CD}
R_{73} @>\sigma_1>> R_{74}\\
 @V\sigma_2VV 	 @VV\sigma_2V \\
R_{75} @>>\sigma_1> R_{76}
\end{CD}
\quad\quad\quad
\begin{aligned}
R_{73} &=	( 0,\ y,\ yz+y+z) \\
R_{74} &=	( xy+x+y,\ y,\ 0) \\
R_{75} &=	( 1,\ y,\ yz) \\
R_{76} &=	( xy,\ y,\ 1) 
\end{aligned}
\end{equation}
\item \begin{equation}
\begin{CD}
R_{77} @>\sigma_1>> R_{78}\\
 @V\sigma_2VV 	 @VV\sigma_2V \\
R_{79} @>>\sigma_1> R_{80}
\end{CD}
\quad\quad\quad
\begin{aligned}
R_{77} &=	( 0,\ z,\ yz) \\
R_{78} &=	( xy,\ x,\ 0) \\
R_{79} &=	( 1,\ z,\ yz+y+z) \\
R_{80} &=	( xy+x+y,\ x,\ 1) 
\end{aligned}
\end{equation}
\item \begin{equation}
\begin{CD}
R_{81} @>\sigma_1>> R_{82}\\
 @V\sigma_2VV 	 @VV\sigma_2V \\
R_{83} @>>\sigma_1> R_{84}
\end{CD}
\quad\quad\quad
\begin{aligned}
R_{81} &=	( 0,\ z,\ yz+y+z) \\
R_{82} &=	( xy+x+y,\ x,\ 0) \\
R_{83} &=	( 1,\ z,\ yz) \\
R_{84} &=	( xy,\ x,\ 1) 
\end{aligned}
\end{equation}
\item \begin{equation}
\begin{CD}
R_{85} @>\sigma_1>> R_{86}\\
 @V\sigma_2VV 	 @VV\sigma_2V \\
R_{87} @>>\sigma_1> R_{88}
\end{CD}
\quad\quad\quad
\begin{aligned}
R_{85} &=	( 0,\ xy+x+y,\ y) \\
R_{86} &=	( y,\ yz+y+z,\ 0) \\
R_{87} &=	( 1,\ xy,\ y) \\
R_{88} &=	( y,\ yz,\ 1) 
\end{aligned}
\end{equation}
\item \begin{equation}
\begin{CD}
R_{89} @>\sigma_1>> R_{90}\\
 @V\sigma_2VV 	 @VV\sigma_2V \\
R_{91} @>>\sigma_1> R_{92}
\end{CD}
\quad\quad\quad
\begin{aligned}
R_{89} &=	( 0,\ xz,\ z) \\
R_{90} &=	( x,\ xz,\ 0) \\
R_{91} &=	( 1,\ xz+x+z,\ z) \\
R_{92} &=	( x,\ xz+x+z,\ 1) 
\end{aligned}
\end{equation}
\item \begin{equation}
\begin{CD}
R_{93} @>\sigma_1>> R_{94}\\
 @V\sigma_2VV 	 @VV\sigma_2V \\
R_{95} @>>\sigma_1> R_{96}
\end{CD}
\quad\quad\quad
\begin{aligned}
R_{93} &=	( 0,\ xz+x+z,\ z) \\
R_{94} &=	( x,\ xz+x+z,\ 0) \\
R_{95} &=	( 1,\ xz,\ z) \\
R_{96} &=	( x,\ xz,\ 1) 
\end{aligned}
\end{equation}
\item \begin{equation}
\begin{CD}
R_{97} @>\sigma_1>> R_{98}\\
 @V\sigma_2VV 	 @VV\sigma_2V \\
R_{99} @>>\sigma_1> R_{100}
\end{CD}
\quad\quad\quad
\begin{aligned}
R_{100} &=	( 0,\ yz,\ y) \\
R_{97} &=	( y,\ xy,\ 0) \\
R_{98} &=	( 1,\ yz+y+z,\ y) \\
R_{99} &=	( y,\ xy+x+y,\ 1) 
\end{aligned}
\end{equation}
\item \begin{equation}
\begin{CD}
R_{101} @>\sigma_1>> R_{102}\\
 @V\sigma_2VV 	 @VV\sigma_2V \\
R_{103} @>>\sigma_1> R_{104}
\end{CD}
\quad\quad\quad
\begin{aligned}
R_{101} &=	( 0,\ yz,\ z) \\
R_{102} &=	( x,\ xy,\ 0) \\
R_{103} &=	( 1,\ yz+y+z,\ z) \\
R_{104} &=	( x,\ xy+x+y,\ 1) 
\end{aligned}
\end{equation}
\item \begin{equation}
\begin{CD}
R_{105} @>\sigma_1>> R_{106}\\
 @V\sigma_2VV 	 @VV\sigma_2V \\
R_{107} @>>\sigma_1> R_{108}
\end{CD}
\quad\quad\quad
\begin{aligned}
R_{105} &=	( 0,\ yz,\ yz+y+z) \\
R_{106} &=	( xy+x+y,\ xy,\ 0) \\
R_{107} &=	( 1,\ yz+y+z,\ yz) \\
R_{108} &=	( xy,\ xy+x+y,\ 1) 
\end{aligned}
\end{equation}
\item \begin{equation}
\begin{CD}
R_{109} @>\sigma_1>> R_{110}\\
 @V\sigma_2VV 	 @VV\sigma_2V \\
R_{111} @>>\sigma_1> R_{112}
\end{CD}
\quad\quad\quad
\begin{aligned}
R_{109} &=	( 0,\ yz+y,\ z) \\
R_{110} &=	( x,\ xy+y,\ 0) \\
R_{111} &=	( 1,\ yz+z+1,\ z) \\
R_{112} &=	( x,\ xy+x+1,\ 1) 
\end{aligned}
\end{equation}
\item \begin{equation}
\begin{CD}
R_{113} @>\sigma_1>> R_{114}\\
 @V\sigma_2VV 	 @VV\sigma_2V \\
R_{115} @>>\sigma_1> R_{116}
\end{CD}
\quad\quad\quad
\begin{aligned}
R_{113} &=	( 0,\ yz+y+z,\ y) \\
R_{114} &=	( y,\ xy+x+y,\ 0) \\
R_{115} &=	( 1,\ yz,\ y) \\
R_{116} &=	( y,\ xy,\ 1) 
\end{aligned}
\end{equation}
\item \begin{equation}
\begin{CD}
R_{117} @>\sigma_1>> R_{118}\\
 @V\sigma_2VV 	 @VV\sigma_2V \\
R_{119} @>>\sigma_1> R_{120}
\end{CD}
\quad\quad\quad
\begin{aligned}
R_{117} &=	( 0,\ yz+y+z,\ z) \\
R_{118} &=	( x,\ xy+x+y,\ 0) \\
R_{119} &=	( 1,\ yz,\ z) \\
R_{120} &=	( x,\ xy,\ 1) 
\end{aligned}
\end{equation}
\item \begin{equation}
\begin{CD}
R_{121} @>\sigma_1>> R_{122}\\
 @V\sigma_2VV 	 @VV\sigma_2V \\
R_{123} @>>\sigma_1> R_{124}
\end{CD}
\quad\quad\quad
\begin{aligned}
R_{121} &=	( 0,\ yz+y+z,\ yz) \\
R_{122} &=	( xy,\ xy+x+y,\ 0) \\
R_{123} &=	( 1,\ yz,\ yz+y+z) \\
R_{124} &=	( xy+x+y,\ xy,\ 1) 
\end{aligned}
\end{equation}
\item \begin{equation}
\begin{CD}
R_{125} @>\sigma_1>> R_{126}\\
 @V\sigma_2VV 	 @VV\sigma_2V \\
R_{127} @>>\sigma_1> R_{128}
\end{CD}
\quad%\quad\quad
\begin{aligned}
R_{125} &=	( 0,\ xyz+xy+xz+yz+x+y+z,\ y) \\
R_{126} &=	( y,\ xyz+xy+xz+yz+x+y+z,\ 0) \\
R_{127} &=	( 1,\ xyz,\ y) \\
R_{128} &=	( y,\ xyz,\ 1) 
\end{aligned}
\end{equation}
\item \begin{equation}
\begin{CD}
R_{129} @>\sigma_1>> R_{130}\\
 @V\sigma_2VV 	 @VV\sigma_2V \\
R_{131} @>>\sigma_1> R_{132}
\end{CD}
\quad%\quad\quad
\begin{aligned}
R_{129} &=	( 0,\ xyz+xy+xz+yz+x+y+z,\ z) \\
R_{130} &=	( x,\ xyz+xy+xz+yz+x+y+z,\ 0) \\
R_{131} &=	( 1,\ xyz,\ z) \\
R_{132} &=	( x,\ xyz,\ 1) 
\end{aligned}
\end{equation}
\item \begin{equation}
\begin{CD}
R_{133} @>\sigma_1>> R_{134}\\
 @V\sigma_2VV 	 @VV\sigma_2V \\
R_{135} @>>\sigma_1> R_{136}
\end{CD}
\quad%\quad\quad
\begin{aligned}
	R_{133} &=\\	( 0,&xyz+xy+xz+yz+x+y+z,yz+y+z)\\
	R_{134} &=\\	( xy&+x+y,xyz+xy+xz+yz+x+y+z,0)\\
R_{135} &=	( 1,\ xyz,\ yz) \\
R_{136} &=	( xy,\ xyz,\ 1) 
\end{aligned}
\end{equation}
\item \begin{equation}
\begin{CD}
R_{137} @>\sigma_1>> R_{138}\\
 @V\sigma_2VV 	 @VV\sigma_2V \\
R_{139} @>>\sigma_1> R_{140}
\end{CD}
\quad\quad\quad
\begin{aligned}
R_{137} &=	( y,\ 0,\ yz+y+z) \\
R_{138} &=	( xy+x+y,\ 0,\ y) \\
R_{139} &=	( y,\ 1,\ yz) \\
R_{140} &=	( xy,\ 1,\ y) 
\end{aligned}
\end{equation}
\item \begin{equation}
\begin{CD}
R_{141} @>\sigma_1>> R_{142}\\
 @V\sigma_2VV 	 @VV\sigma_2V \\
R_{143} @>>\sigma_1> R_{144}
\end{CD}
\quad\quad\quad
\begin{aligned}
R_{141} &=	( y,\ y,\ yz) \\
R_{142} &=	( xy,\ y,\ y) \\
R_{143} &=	( y,\ y,\ yz+y+z) \\
R_{144} &=	( xy+x+y,\ y,\ y) 
\end{aligned}
\end{equation}
\item \begin{equation}
\begin{CD}
R_{145} @>\sigma_1>> R_{146}\\
 @V\sigma_2VV 	 @VV\sigma_2V \\
R_{147} @>>\sigma_1> R_{148}
\end{CD}
\quad\quad\quad
\begin{aligned}
R_{145} &=	( y,\ y,\ yz+y+1) \\
R_{146} &=	( xy+y+1,\ y,\ y) \\
R_{147} &=	( y,\ y,\ yz+z) \\
R_{148} &=	( xy+x,\ y,\ y) 
\end{aligned}
\end{equation}
\item \begin{equation}
\begin{CD}
R_{149} @>\sigma_1>> R_{150}\\
 @V\sigma_2VV 	 @VV\sigma_2V \\
R_{151} @>>\sigma_1> R_{152}
\end{CD}
\quad\quad\quad
\begin{aligned}
R_{149} &=	( y,\ xy,\ y) \\
R_{150} &=	( y,\ yz,\ y) \\
R_{151} &=	( y,\ xy+x+y,\ y) \\
R_{152} &=	( y,\ yz+y+z,\ y) 
\end{aligned}
\end{equation}
\item \begin{equation}
\begin{CD}
R_{153} @>\sigma_1>> R_{154}\\
 @V\sigma_2VV 	 @VV\sigma_2V \\
R_{155} @>>\sigma_1> R_{156}
\end{CD}
\quad\quad\quad
\begin{aligned}
R_{153} &=	( y,\ yz,\ yz) \\
R_{154} &=	( xy,\ xy,\ y) \\
R_{155} &=	( y,\ yz+y+z,\ yz+y+z) \\
R_{156} &=	( xy+x+y,\ xy+x+y,\ y) 
\end{aligned}
\end{equation}
\item \begin{equation}
\begin{CD}
R_{157}\\
 @V\sigma_2VV \\
R_{158}
\end{CD}
\quad\quad\quad
\begin{aligned}
R_{157} &=	( y,\ xyz,\ y) \\
R_{158} &=	( y,\ xyz+xy+xz+yz+x+y+z,\ y) 
\end{aligned}
\end{equation}
\item \begin{equation}
\begin{CD}
R_{159} @>\sigma_1>> R_{160}\\
 @V\sigma_2VV 	 @VV\sigma_2V \\
R_{161} @>>\sigma_1> R_{162}
\end{CD}
\quad\quad\quad
\begin{aligned}
R_{159} &=	( xy,\ xy,\ xy+y) \\
R_{160} &=	( yz+y,\ yz,\ yz) \\
R_{161} &=	( xy+x+y,\ xy+x+y,\ xy+x+1) \\
R_{162} &=	( yz+z+1,\ yz+y+z,\ yz+y+z) 
\end{aligned}
\end{equation}

\subsection{Solutions of image cardinality 4}
\item \begin{equation}
\begin{CD}
R_{163} @>\sigma_1>> R_{164}\\
 @V\sigma_2VV 	 @VV\sigma_2V \\
R_{165} @>>\sigma_1> R_{166}
\end{CD}
\quad\quad\quad
\begin{aligned}
R_{163} &=	( 0,\ x,\ y) \\
R_{164} &=	( y,\ z,\ 0) \\
R_{165} &=	( 1,\ x,\ y) \\
R_{166} &=	( y,\ z,\ 1) 
\end{aligned}
\end{equation}
\item \begin{equation}
\begin{CD}
R_{167} @>\sigma_1>> R_{168}\\
 @V\sigma_2VV 	 @VV\sigma_2V \\
R_{169} @>>\sigma_1> R_{170}
\end{CD}
\quad\quad\quad
\begin{aligned}
R_{167} &=	( 0,\ x,\ z) \\
R_{168} &=	( x,\ z,\ 0) \\
R_{169} &=	( 1,\ x,\ z) \\
R_{170} &=	( x,\ z,\ 1) 
\end{aligned}
\end{equation}
\item \begin{equation}
\begin{CD}
R_{171} @>\sigma_1>> R_{172}\\
 @V\sigma_2VV 	 @VV\sigma_2V \\
R_{173} @>>\sigma_1> R_{174}
\end{CD}
\quad\quad\quad
\begin{aligned}
R_{171} &=	( 0,\ x,\ y+z) \\
R_{172} &=	( x+y,\ z,\ 0) \\
R_{173} &=	( 1,\ x,\ y+z+1) \\
R_{174} &=	( x+y+1,\ z,\ 1) 
\end{aligned}
\end{equation}
\begin{multline}
{\bf R_{171}\colon}\quad \im\delta^2=
	\begin{pmatrix}
0 & 0 & 0 & -2 & 0 & 0 & 0 & -2 
	\end{pmatrix}^{\rm T}
\\ \ker\delta^3=
	\begin{pmatrix}
1 & 1 & 1 & 0 & 1 & 1 & 1 & 0 \\
0 & 0 & 0 & 1 & 0 & 0 & 0 & 1 
	\end{pmatrix}^{\rm T}
\end{multline}
\begin{multline}
{\bf R_{172}\colon}\quad \im\delta^2=
	\begin{pmatrix}
0 & 0 & 0 & 0 & 0 & 0 & -2 & -2 
	\end{pmatrix}^{\rm T}
\\ \ker\delta^3=
	\begin{pmatrix}
1 & 1 & 1 & 1 & 1 & 1 & 0 & 0 \\
0 & 0 & 0 & 0 & 0 & 0 & 1 & 1 
	\end{pmatrix}^{\rm T}
\end{multline}
\begin{multline}
{\bf R_{173}\colon}\quad \im\delta^2=
	\begin{pmatrix}
2 & 0 & 0 & 0 & 2 & 0 & 0 & 0 
	\end{pmatrix}^{\rm T}
\\ \ker\delta^3=
	\begin{pmatrix}
1 & 0 & 0 & 0 & 1 & 0 & 0 & 0 \\
0 & 1 & 1 & 1 & 0 & 1 & 1 & 1 
	\end{pmatrix}^{\rm T}
\end{multline}
\begin{multline}
{\bf R_{174}\colon}\quad \im\delta^2=
	\begin{pmatrix}
2 & 2 & 0 & 0 & 0 & 0 & 0 & 0 
	\end{pmatrix}^{\rm T}
\\ \ker\delta^3=
	\begin{pmatrix}
1 & 1 & 0 & 0 & 0 & 0 & 0 & 0 \\
0 & 0 & 1 & 1 & 1 & 1 & 1 & 1 
	\end{pmatrix}^{\rm T}
\end{multline}
\item \begin{equation}
\begin{CD}
R_{175} @>\sigma_1>> R_{176}\\
 @V\sigma_2VV 	 @VV\sigma_2V \\
R_{177} @>>\sigma_1> R_{178}
\end{CD}
\quad\quad\quad
\begin{aligned}
R_{175} &=	( 0,\ x,\ x+y+z) \\
R_{176} &=	( x+y+z,\ z,\ 0) \\
R_{177} &=	( 1,\ x,\ x+y+z) \\
R_{178} &=	( x+y+z,\ z,\ 1) 
\end{aligned}
\end{equation}
\item \begin{equation}
\begin{CD}
R_{179} @>\sigma_1>> R_{180}\\
 @V\sigma_2VV 	 @VV\sigma_2V \\
R_{181} @>>\sigma_1> R_{182}
\end{CD}
\quad\quad\quad
\begin{aligned}
R_{179} &=	( 0,\ x,\ yz+y+z) \\
R_{180} &=	( xy+x+y,\ z,\ 0) \\
R_{181} &=	( 1,\ x,\ yz) \\
R_{182} &=	( xy,\ z,\ 1) 
\end{aligned}
\end{equation}
\item \begin{equation}
\begin{CD}
R_{183} @>\sigma_1>> R_{184}\\
 @V\sigma_2VV 	 @VV\sigma_2V \\
R_{185} @>>\sigma_1> R_{186}
\end{CD}
\quad\quad\quad
\begin{aligned}
R_{183} &=	( 0,\ y,\ z) \\
R_{184} &=	( x,\ y,\ 0) \\
R_{185} &=	( 1,\ y,\ z) \\
R_{186} &=	( x,\ y,\ 1) 
\end{aligned}
\end{equation}
\item \begin{equation}
\begin{CD}
R_{187} @>\sigma_1>> R_{188}\\
 @V\sigma_2VV 	 @VV\sigma_2V \\
R_{189} @>>\sigma_1> R_{190}
\end{CD}
\quad\quad\quad
\begin{aligned}
R_{187} &=	( 0,\ x+y,\ y) \\
R_{188} &=	( y,\ y+z,\ 0) \\
R_{189} &=	( 1,\ x+y+1,\ y) \\
R_{190} &=	( y,\ y+z+1,\ 1) 
\end{aligned}
\end{equation}
\item \begin{equation}
\begin{CD}
R_{191} @>\sigma_1>> R_{192}\\
 @V\sigma_2VV 	 @VV\sigma_2V \\
R_{193} @>>\sigma_1> R_{194}
\end{CD}
\quad\quad\quad
\begin{aligned}
R_{191} &=	( 0,\ x+y,\ z) \\
R_{192} &=	( x,\ y+z,\ 0) \\
R_{193} &=	( 1,\ x+y+1,\ z) \\
R_{194} &=	( x,\ y+z+1,\ 1) 
\end{aligned}
\end{equation}
\begin{multline}
{\bf R_{191}\colon}\quad \im\delta^2=
	\begin{pmatrix}
0 & 0 & 0 & 0 & 0 & 0 & -2 & -2 
	\end{pmatrix}^{\rm T}
\\ \ker\delta^3=
	\begin{pmatrix}
1 & 1 & 1 & 1 & 1 & 1 & 0 & 0 \\
0 & 0 & 0 & 0 & 0 & 0 & 1 & 1 
	\end{pmatrix}^{\rm T}
\end{multline}
\begin{multline}
{\bf R_{192}\colon}\quad \im\delta^2=
	\begin{pmatrix}
0 & 0 & 0 & -2 & 0 & 0 & 0 & -2 
	\end{pmatrix}^{\rm T}
\\ \ker\delta^3=
	\begin{pmatrix}
1 & 1 & 1 & 0 & 1 & 1 & 1 & 0 \\
0 & 0 & 0 & 1 & 0 & 0 & 0 & 1 
	\end{pmatrix}^{\rm T}
\end{multline}
\begin{multline}
{\bf R_{193}\colon}\quad \im\delta^2=
	\begin{pmatrix}
2 & 2 & 0 & 0 & 0 & 0 & 0 & 0 
	\end{pmatrix}^{\rm T}
\\ \ker\delta^3=
	\begin{pmatrix}
1 & 1 & 0 & 0 & 0 & 0 & 0 & 0 \\
0 & 0 & 1 & 1 & 1 & 1 & 1 & 1 
	\end{pmatrix}^{\rm T}
\end{multline}
\begin{multline}
{\bf R_{194}\colon}\quad \im\delta^2=
	\begin{pmatrix}
2 & 0 & 0 & 0 & 2 & 0 & 0 & 0 
	\end{pmatrix}^{\rm T}
\\ \ker\delta^3=
	\begin{pmatrix}
1 & 0 & 0 & 0 & 1 & 0 & 0 & 0 \\
0 & 1 & 1 & 1 & 0 & 1 & 1 & 1 
	\end{pmatrix}^{\rm T}
\end{multline}
\item \begin{equation}
\begin{CD}
R_{195} @>\sigma_1>> R_{196}\\
 @V\sigma_2VV 	 @VV\sigma_2V \\
R_{197} @>>\sigma_1> R_{198}
\end{CD}
\quad\quad\quad
\begin{aligned}
R_{195} &=	( 0,\ x+y,\ x+z) \\
R_{196} &=	( x+z,\ y+z,\ 0) \\
R_{197} &=	( 1,\ x+y+1,\ x+z+1) \\
R_{198} &=	( x+z+1,\ y+z+1,\ 1) 
\end{aligned}
\end{equation}
\item \begin{equation}
\begin{CD}
R_{199} @>\sigma_1>> R_{200}\\
 @V\sigma_2VV 	 @VV\sigma_2V \\
R_{201} @>>\sigma_1> R_{202}
\end{CD}
\quad\quad\quad
\begin{aligned}
R_{199} &=	( 0,\ z,\ y) \\
R_{200} &=	( y,\ x,\ 0) \\
R_{201} &=	( 1,\ z,\ y) \\
R_{202} &=	( y,\ x,\ 1) 
\end{aligned}
\end{equation}
\item \begin{equation}
\begin{CD}
R_{203} @>\sigma_1>> R_{204}\\
 @V\sigma_2VV 	 @VV\sigma_2V \\
R_{205} @>>\sigma_1> R_{206}
\end{CD}
\quad\quad\quad
\begin{aligned}
R_{203} &=	( 0,\ x+z,\ y) \\
R_{204} &=	( y,\ x+z,\ 0) \\
R_{205} &=	( 1,\ x+z+1,\ y) \\
R_{206} &=	( y,\ x+z+1,\ 1) 
\end{aligned}
\end{equation}
\begin{multline}
{\bf R_{203}\colon}\quad \im\delta^2=
	\begin{pmatrix}
0 & 0 & 0 & 0 & 0 & -2 & 0 & -2 
	\end{pmatrix}^{\rm T}
\\ \ker\delta^3=
	\begin{pmatrix}
1 & 1 & 1 & 1 & 1 & 0 & 1 & 0 \\
0 & 0 & 0 & 0 & 0 & 1 & 0 & 1 
	\end{pmatrix}^{\rm T}
\end{multline}
\begin{multline}
{\bf R_{204}\colon}\quad \im\delta^2=
	\begin{pmatrix}
0 & 0 & 0 & 0 & 0 & -2 & 0 & -2 
	\end{pmatrix}^{\rm T}
\\ \ker\delta^3=
	\begin{pmatrix}
1 & 1 & 1 & 1 & 1 & 0 & 1 & 0 \\
0 & 0 & 0 & 0 & 0 & 1 & 0 & 1 
	\end{pmatrix}^{\rm T}
\end{multline}
\begin{multline}
{\bf R_{205}\colon}\quad \im\delta^2=
	\begin{pmatrix}
2 & 0 & 2 & 0 & 0 & 0 & 0 & 0 
	\end{pmatrix}^{\rm T}
\\ \ker\delta^3=
	\begin{pmatrix}
1 & 0 & 1 & 0 & 0 & 0 & 0 & 0 \\
0 & 1 & 0 & 1 & 1 & 1 & 1 & 1 
	\end{pmatrix}^{\rm T}
\end{multline}
\begin{multline}
{\bf R_{206}\colon}\quad \im\delta^2=
	\begin{pmatrix}
2 & 0 & 2 & 0 & 0 & 0 & 0 & 0 
	\end{pmatrix}^{\rm T}
\\ \ker\delta^3=
	\begin{pmatrix}
1 & 0 & 1 & 0 & 0 & 0 & 0 & 0 \\
0 & 1 & 0 & 1 & 1 & 1 & 1 & 1 
	\end{pmatrix}^{\rm T}
\end{multline}
\item \begin{equation}
\begin{CD}
R_{207} @>\sigma_1>> R_{208}\\
 @V\sigma_2VV 	 @VV\sigma_2V \\
R_{209} @>>\sigma_1> R_{210}
\end{CD}
\quad\quad\quad
\begin{aligned}
R_{207} &=	( 0,\ x+z,\ z) \\
R_{208} &=	( x,\ x+z,\ 0) \\
R_{209} &=	( 1,\ x+z+1,\ z) \\
R_{210} &=	( x,\ x+z+1,\ 1) 
\end{aligned}
\end{equation}
\item \begin{equation}
\begin{CD}
R_{211} @>\sigma_1>> R_{212}\\
 @V\sigma_2VV 	 @VV\sigma_2V \\
R_{213} @>>\sigma_1> R_{214}
\end{CD}
\quad\quad\quad
\begin{aligned}
R_{211} &=	( 0,\ xy+x+y,\ z) \\
R_{212} &=	( x,\ yz+y+z,\ 0) \\
R_{213} &=	( 1,\ xy,\ z) \\
R_{214} &=	( x,\ yz,\ 1) 
\end{aligned}
\end{equation}
\item \begin{equation}
\begin{CD}
R_{215} @>\sigma_1>> R_{216}\\
 @V\sigma_2VV 	 @VV\sigma_2V \\
R_{217} @>>\sigma_1> R_{218}
\end{CD}
\quad\quad\quad
\begin{aligned}
R_{215} &=	( 0,\ xy+x+y,\ yz+y+z) \\
R_{216} &=	( xy+x+y,\ yz+y+z,\ 0) \\
R_{217} &=	( 1,\ xy,\ yz) \\
R_{218} &=	( xy,\ yz,\ 1) 
\end{aligned}
\end{equation}
\item \begin{equation}
\begin{CD}
R_{219} @>\sigma_1>> R_{220}\\
 @V\sigma_2VV 	 @VV\sigma_2V \\
R_{221} @>>\sigma_1> R_{222}
\end{CD}
\quad\quad\quad
\begin{aligned}
R_{219} &=	( 0,\ xz+x+z,\ y) \\
R_{220} &=	( y,\ xz+x+z,\ 0) \\
R_{221} &=	( 1,\ xz,\ y) \\
R_{222} &=	( y,\ xz,\ 1) 
\end{aligned}
\end{equation}
\item \begin{equation}
\begin{CD}
R_{223} @>\sigma_1>> R_{224}\\
 @V\sigma_2VV 	 @VV\sigma_2V \\
R_{225} @>>\sigma_1> R_{226}
\end{CD}
\quad\quad\quad
\begin{aligned}
R_{223} &=	( 0,\ xz+x+z,\ yz+y+z) \\
R_{224} &=	( xy+x+y,\ xz+x+z,\ 0) \\
R_{225} &=	( 1,\ xz,\ yz) \\
R_{226} &=	( xy,\ xz,\ 1) 
\end{aligned}
\end{equation}
\item \begin{equation}
\begin{CD}
R_{227} @>\sigma_1>> R_{228}\\
 @V\sigma_2VV 	 @VV\sigma_2V \\
R_{229} @>>\sigma_1> R_{230}
\end{CD}
\quad\quad\quad
\begin{aligned}
R_{227} &=	( x,\ 0,\ y) \\
R_{228} &=	( y,\ 0,\ z) \\
R_{229} &=	( x,\ 1,\ y) \\
R_{230} &=	( y,\ 1,\ z) 
\end{aligned}
\end{equation}
\item \begin{equation}
\begin{CD}
R_{231}\\
 @V\sigma_2VV \\
R_{232}
\end{CD}
\quad\quad\quad
\begin{aligned}
R_{231} &=	( x,\ 0,\ z) \\
R_{232} &=	( x,\ 1,\ z) 
\end{aligned}
\end{equation}
\item \begin{equation}
\begin{CD}
R_{233} @>\sigma_1>> R_{234}\\
 @V\sigma_2VV 	 @VV\sigma_2V \\
R_{235} @>>\sigma_1> R_{236}
\end{CD}
\quad\quad\quad
\begin{aligned}
R_{233} &=	( x,\ 0,\ y+z) \\
R_{234} &=	( x+y,\ 0,\ z) \\
R_{235} &=	( x,\ 1,\ y+z+1) \\
R_{236} &=	( x+y+1,\ 1,\ z) 
\end{aligned}
\end{equation}
\begin{multline}
{\bf R_{233}\colon}\quad \im\delta^2=
	\begin{pmatrix}
0 & 0 & 0 & -2 & 0 & 0 & 0 & -2 
	\end{pmatrix}^{\rm T}
\\ \ker\delta^3=
	\begin{pmatrix}
1 & 1 & 1 & 0 & 1 & 1 & 1 & 0 \\
0 & 0 & 0 & 1 & 0 & 0 & 0 & 1 
	\end{pmatrix}^{\rm T}
\end{multline}
\begin{multline}
{\bf R_{234}\colon}\quad \im\delta^2=
	\begin{pmatrix}
0 & 0 & 0 & 0 & 0 & 0 & -2 & -2 
	\end{pmatrix}^{\rm T}
\\ \ker\delta^3=
	\begin{pmatrix}
1 & 1 & 1 & 1 & 1 & 1 & 0 & 0 \\
0 & 0 & 0 & 0 & 0 & 0 & 1 & 1 
	\end{pmatrix}^{\rm T}
\end{multline}
\begin{multline}
{\bf R_{235}\colon}\quad \im\delta^2=
	\begin{pmatrix}
2 & 0 & 0 & 0 & 2 & 0 & 0 & 0 
	\end{pmatrix}^{\rm T}
\\ \ker\delta^3=
	\begin{pmatrix}
1 & 0 & 0 & 0 & 1 & 0 & 0 & 0 \\
0 & 1 & 1 & 1 & 0 & 1 & 1 & 1 
	\end{pmatrix}^{\rm T}
\end{multline}
\begin{multline}
{\bf R_{236}\colon}\quad \im\delta^2=
	\begin{pmatrix}
2 & 2 & 0 & 0 & 0 & 0 & 0 & 0 
	\end{pmatrix}^{\rm T}
\\ \ker\delta^3=
	\begin{pmatrix}
1 & 1 & 0 & 0 & 0 & 0 & 0 & 0 \\
0 & 0 & 1 & 1 & 1 & 1 & 1 & 1 
	\end{pmatrix}^{\rm T}
\end{multline}
\item \begin{equation}
\begin{CD}
R_{237} @>\sigma_1>> R_{238}\\
 @V\sigma_2VV 	 @VV\sigma_2V \\
R_{239} @>>\sigma_1> R_{240}
\end{CD}
\quad\quad\quad
\begin{aligned}
R_{237} &=	( x,\ 0,\ yz+y+z) \\
R_{238} &=	( xy+x+y,\ 0,\ z) \\
R_{239} &=	( x,\ 1,\ yz) \\
R_{240} &=	( xy,\ 1,\ z) 
\end{aligned}
\end{equation}
\item \begin{equation}
\begin{CD}
R_{241} @>\sigma_1>> R_{242}
\end{CD}
\quad\quad\quad
\begin{aligned}
R_{241} &=	( x,\ x,\ y) \\
R_{242} &=	( y,\ z,\ z) 
\end{aligned}
\end{equation}
\item \begin{equation}
\begin{CD}
R_{243} @>\sigma_1>> R_{244}\\
 @V\sigma_2VV 	 @VV\sigma_2V \\
R_{245} @>>\sigma_1> R_{246}
\end{CD}
\quad\quad\quad
\begin{aligned}
R_{243} &=	( x,\ x,\ x+y) \\
R_{244} &=	( y+z,\ z,\ z) \\
R_{245} &=	( x,\ x,\ x+y+1) \\
R_{246} &=	( y+z+1,\ z,\ z) 
\end{aligned}
\end{equation}
\item \begin{equation}
\begin{CD}
R_{247} @>\sigma_1>> R_{248}
\end{CD}
\quad\quad\quad
\begin{aligned}
R_{247} &=	( x,\ x,\ z) \\
R_{248} &=	( x,\ z,\ z) 
\end{aligned}
\end{equation}
\item \begin{equation}
\begin{CD}
R_{249} @>\sigma_1>> R_{250}
\end{CD}
\quad\quad\quad
\begin{aligned}
R_{249} &=	( x,\ x,\ z+1) \\
R_{250} &=	( x+1,\ z,\ z) 
\end{aligned}
\end{equation}
\item \begin{equation}
\begin{CD}
R_{251} @>\sigma_1>> R_{252}
\end{CD}
\quad\quad\quad
\begin{aligned}
R_{251} &=	( x,\ x,\ x+y+z) \\
R_{252} &=	( x+y+z,\ z,\ z) 
\end{aligned}
\end{equation}
\begin{multline}
{\bf R_{251}\colon}\quad \im\delta^2=
	\begin{pmatrix}
0 & 0 & 0 & -2 & 2 & 0 & 0 & 0 
	\end{pmatrix}^{\rm T}
\\ \ker\delta^3=
	\begin{pmatrix}
1 & 1 & 1 & 0 & 2 & 1 & 1 & 1 \\
0 & 0 & 0 & 1 & -1 & 0 & 0 & 0 
	\end{pmatrix}^{\rm T}
\end{multline}
\begin{multline}
{\bf R_{252}\colon}\quad \im\delta^2=
	\begin{pmatrix}
0 & 2 & 0 & 0 & 0 & 0 & -2 & 0 
	\end{pmatrix}^{\rm T}
\\ \ker\delta^3=
	\begin{pmatrix}
1 & 0 & 1 & 1 & 1 & 1 & 2 & 1 \\
0 & 1 & 0 & 0 & 0 & 0 & -1 & 0 
	\end{pmatrix}^{\rm T}
\end{multline}
\item \begin{equation}
\begin{CD}
R_{253} @>\sigma_1>> R_{254}
\end{CD}
\quad\quad\quad
\begin{aligned}
R_{253} &=	( x,\ x,\ x+y+z+1) \\
R_{254} &=	( x+y+z+1,\ z,\ z) 
\end{aligned}
\end{equation}
\begin{multline}
{\bf R_{253}\colon}\quad \im\delta^2=
	\begin{pmatrix}
1 & -1 & -1 & -1 & 1 & 1 & 1 & -1 
	\end{pmatrix}^{\rm T}
\\ \ker\delta^3=
	\begin{pmatrix}
1 & 0 & 0 & 0 & 1 & 1 & 1 & 0 \\
0 & 1 & 1 & 1 & 0 & 0 & 0 & 1 
	\end{pmatrix}^{\rm T}
\end{multline}
\begin{multline}
{\bf R_{254}\colon}\quad \im\delta^2=
	\begin{pmatrix}
1 & 1 & -1 & 1 & -1 & 1 & -1 & -1 
	\end{pmatrix}^{\rm T}
\\ \ker\delta^3=
	\begin{pmatrix}
1 & 1 & 0 & 1 & 0 & 1 & 0 & 0 \\
0 & 0 & 1 & 0 & 1 & 0 & 1 & 1 
	\end{pmatrix}^{\rm T}
\end{multline}
\item \begin{equation}
\begin{CD}
R_{255} @>\sigma_1>> R_{256}
\end{CD}
\quad\quad\quad
\begin{aligned}
R_{255} &=	( x,\ y,\ y) \\
R_{256} &=	( y,\ y,\ z) 
\end{aligned}
\end{equation}
\item \begin{equation}
\begin{CD}
R_{257} @>\sigma_1>> R_{258}\\
 @V\sigma_2VV 	 @VV\sigma_2V \\
R_{259} @>>\sigma_1> R_{260}
\end{CD}
\quad\quad\quad
\begin{aligned}
R_{257} &=	( x,\ y,\ x+y) \\
R_{258} &=	( y+z,\ y,\ z) \\
R_{259} &=	( x,\ y,\ x+y+1) \\
R_{260} &=	( y+z+1,\ y,\ z) 
\end{aligned}
\end{equation}
\item \begin{equation}
\begin{CD}
R_{261} @>\sigma_1>> R_{262}
\end{CD}
\quad\quad\quad
\begin{aligned}
R_{261} &=	( x,\ y+1,\ y) \\
R_{262} &=	( y,\ y+1,\ z) 
\end{aligned}
\end{equation}
\item \begin{equation}
\begin{CD}
R_{263} @>\sigma_1>> R_{264}\\
 @V\sigma_2VV 	 @VV\sigma_2V \\
R_{265} @>>\sigma_1> R_{266}
\end{CD}
\quad\quad\quad
\begin{aligned}
R_{263} &=	( x,\ y+1,\ x+y) \\
R_{264} &=	( y+z,\ y+1,\ z) \\
R_{265} &=	( x,\ y+1,\ x+y+1) \\
R_{266} &=	( y+z+1,\ y+1,\ z) 
\end{aligned}
\end{equation}
\item \begin{equation}
\begin{CD}
R_{267} @>\sigma_1>> R_{268}
\end{CD}
\quad\quad\quad
\begin{aligned}
R_{267} &=	( x,\ z,\ z+1) \\
R_{268} &=	( x+1,\ x,\ z) 
\end{aligned}
\end{equation}
\item \begin{equation}
\begin{CD}
R_{269}\\
 @V\sigma_2VV \\
R_{270}
\end{CD}
\quad\quad\quad
\begin{aligned}
R_{269} &=	( x,\ x+z,\ z) \\
R_{270} &=	( x,\ x+z+1,\ z) 
\end{aligned}
\end{equation}
\item \begin{equation}
\begin{CD}
R_{271} @>\sigma_1>> R_{272}\\
 @V\sigma_2VV 	 @VV\sigma_2V \\
R_{273} @>>\sigma_1> R_{274}
\end{CD}
\quad\quad\quad
\begin{aligned}
R_{271} &=	( x,\ xy,\ y) \\
R_{272} &=	( y,\ yz,\ z) \\
R_{273} &=	( x,\ xy+x+y,\ y) \\
R_{274} &=	( y,\ yz+y+z,\ z) 
\end{aligned}
\end{equation}
\item \begin{equation}
\begin{CD}
R_{275} @>\sigma_1>> R_{276}\\
 @V\sigma_2VV 	 @VV\sigma_2V \\
R_{277} @>>\sigma_1> R_{278}
\end{CD}
\quad\quad\quad
\begin{aligned}
R_{275} &=	( x,\ xy+x+1,\ y) \\
R_{276} &=	( y,\ yz+z+1,\ z) \\
R_{277} &=	( x,\ xy+y,\ y) \\
R_{278} &=	( y,\ yz+y,\ z) 
\end{aligned}
\end{equation}
\item \begin{equation}
\begin{CD}
R_{279}\\
 @V\sigma_2VV \\
R_{280}
\end{CD}
\quad\quad\quad
\begin{aligned}
R_{279} &=	( x,\ xz,\ z) \\
R_{280} &=	( x,\ xz+x+z,\ z) 
\end{aligned}
\end{equation}
\item \begin{equation}
\begin{CD}
R_{281} @>\sigma_1>> R_{282}\\
 @V\sigma_2VV 	 @VV\sigma_2V \\
R_{283} @>>\sigma_1> R_{284}
\end{CD}
\quad\quad\quad
\begin{aligned}
R_{281} &=	( x,\ yz,\ yz) \\
R_{282} &=	( xy,\ xy,\ z) \\
R_{283} &=	( x,\ yz+y+z,\ yz+y+z) \\
R_{284} &=	( xy+x+y,\ xy+x+y,\ z) 
\end{aligned}
\end{equation}
\item \begin{equation}
\begin{CD}
R_{285} @>\sigma_1>> R_{286}
\end{CD}
\quad\quad\quad
\begin{aligned}
R_{285} &=	( x+1,\ x,\ y) \\
R_{286} &=	( y,\ z,\ z+1) 
\end{aligned}
\end{equation}
\item \begin{equation}
\begin{CD}
R_{287} @>\sigma_1>> R_{288}
\end{CD}
\quad\quad\quad
\begin{aligned}
R_{287} &=	( x+1,\ x,\ z+1) \\
R_{288} &=	( x+1,\ z,\ z+1) 
\end{aligned}
\end{equation}
\item \begin{equation}
\begin{CD}
R_{289} @>\sigma_1>> R_{290}
\end{CD}
\quad\quad\quad
\begin{aligned}
R_{289} &=	( x+1,\ x,\ x+y+z) \\
R_{290} &=	( x+y+z,\ z,\ z+1) 
\end{aligned}
\end{equation}
\begin{multline}
{\bf R_{289}\colon}\quad \im\delta^2=
	\begin{pmatrix}
1 & 1 & 1 & -1 & 1 & -1 & -1 & -1 
	\end{pmatrix}^{\rm T}
\\ \ker\delta^3=
	\begin{pmatrix}
1 & 1 & 1 & 0 & 1 & 0 & 0 & 0 \\
0 & 0 & 0 & 1 & 0 & 1 & 1 & 1 
	\end{pmatrix}^{\rm T}
\end{multline}
\begin{multline}
{\bf R_{290}\colon}\quad \im\delta^2=
	\begin{pmatrix}
1 & 1 & 1 & -1 & 1 & -1 & -1 & -1 
	\end{pmatrix}^{\rm T}
\\ \ker\delta^3=
	\begin{pmatrix}
1 & 1 & 1 & 0 & 1 & 0 & 0 & 0 \\
0 & 0 & 0 & 1 & 0 & 1 & 1 & 1 
	\end{pmatrix}^{\rm T}
\end{multline}
\item \begin{equation}
\begin{CD}
R_{291} @>\sigma_1>> R_{292}
\end{CD}
\quad\quad\quad
\begin{aligned}
R_{291} &=	( x+1,\ x,\ x+y+z+1) \\
R_{292} &=	( x+y+z+1,\ z,\ z+1) 
\end{aligned}
\end{equation}
\begin{multline}
{\bf R_{291}\colon}\quad \im\delta^2=
	\begin{pmatrix}
2 & 0 & 0 & 0 & 0 & 0 & 0 & -2 
	\end{pmatrix}^{\rm T}
\\ \ker\delta^3=
	\begin{pmatrix}
1 & 0 & 0 & 0 & 0 & 0 & 0 & -1 \\
0 & 1 & 1 & 1 & 1 & 1 & 1 & 2 
	\end{pmatrix}^{\rm T}
\end{multline}
\begin{multline}
{\bf R_{292}\colon}\quad \im\delta^2=
	\begin{pmatrix}
2 & 0 & 0 & 0 & 0 & 0 & 0 & -2 
	\end{pmatrix}^{\rm T}
\\ \ker\delta^3=
	\begin{pmatrix}
1 & 0 & 0 & 0 & 0 & 0 & 0 & -1 \\
0 & 1 & 1 & 1 & 1 & 1 & 1 & 2 
	\end{pmatrix}^{\rm T}
\end{multline}
\item \begin{equation}
\begin{CD}
R_{293} @>\sigma_1>> R_{294}
\end{CD}
\quad\quad\quad
\begin{aligned}
R_{293} &=	( x+1,\ y,\ y) \\
R_{294} &=	( y,\ y,\ z+1) 
\end{aligned}
\end{equation}
\item \begin{equation}
\begin{CD}
R_{295} @>\sigma_1>> R_{296}
\end{CD}
\quad\quad\quad
\begin{aligned}
R_{295} &=	( x+1,\ y+1,\ y) \\
R_{296} &=	( y,\ y+1,\ z+1) 
\end{aligned}
\end{equation}
\item \begin{equation}
\begin{CD}
R_{297} @>\sigma_1>> R_{298}\\
 @V\sigma_2VV 	 @VV\sigma_2V \\
R_{299} @>>\sigma_1> R_{300}
\end{CD}
\quad\quad\quad
\begin{aligned}
R_{297} &=	( y,\ 0,\ y+z) \\
R_{298} &=	( x+y,\ 0,\ y) \\
R_{299} &=	( y,\ 1,\ y+z+1) \\
R_{300} &=	( x+y+1,\ 1,\ y) 
\end{aligned}
\end{equation}
\item \begin{equation}
\begin{CD}
R_{301} @>\sigma_1>> R_{302}
\end{CD}
\quad\quad\quad
\begin{aligned}
R_{301} &=	( y,\ x,\ y) \\
R_{302} &=	( y,\ z,\ y) 
\end{aligned}
\end{equation}
\item \begin{equation}
\begin{CD}
R_{303} @>\sigma_1>> R_{304}\\
 @V\sigma_2VV 	 @VV\sigma_2V \\
R_{305} @>>\sigma_1> R_{306}
\end{CD}
\quad\quad\quad
\begin{aligned}
R_{303} &=	( y,\ x,\ xy+x+1) \\
R_{304} &=	( yz+z+1,\ z,\ y) \\
R_{305} &=	( y,\ x,\ xy+y) \\
R_{306} &=	( yz+y,\ z,\ y) 
\end{aligned}
\end{equation}
\item \begin{equation}
\begin{CD}
R_{307} @>\sigma_1>> R_{308}\\
 @V\sigma_2VV 	 @VV\sigma_2V \\
R_{309} @>>\sigma_1> R_{310}
\end{CD}
\quad\quad\quad
\begin{aligned}
R_{307} &=	( y,\ z,\ yz) \\
R_{308} &=	( xy,\ x,\ y) \\
R_{309} &=	( y,\ z,\ yz+y+z) \\
R_{310} &=	( xy+x+y,\ x,\ y) 
\end{aligned}
\end{equation}
\item \begin{equation}
\begin{CD}
R_{311}\\
 @V\sigma_2VV \\
R_{312}
\end{CD}
\quad\quad\quad
\begin{aligned}
R_{311} &=	( y,\ xz,\ y) \\
R_{312} &=	( y,\ xz+x+z,\ y) 
\end{aligned}
\end{equation}
\item \begin{equation}
\begin{CD}
R_{313} @>\sigma_1>> R_{314}\\
 @V\sigma_2VV 	 @VV\sigma_2V \\
R_{315} @>>\sigma_1> R_{316}
\end{CD}
\quad\quad\quad
\begin{aligned}
R_{313} &=	( y,\ yz,\ yz+y+z) \\
R_{314} &=	( xy+x+y,\ xy,\ y) \\
R_{315} &=	( y,\ yz+y+z,\ yz) \\
R_{316} &=	( xy,\ xy+x+y,\ y) 
\end{aligned}
\end{equation}
\item \begin{equation}
\begin{CD}
R_{317}\\
 @V\sigma_2VV \\
R_{318}
\end{CD}
\quad\quad\quad
\begin{aligned}
R_{317} &=	( x+y,\ 0,\ y+z) \\
R_{318} &=	( x+y+1,\ 1,\ y+z+1) 
\end{aligned}
\end{equation}
\item \begin{equation}
\begin{CD}
R_{319}\\
 @V\sigma_2VV \\
R_{320}
\end{CD}
\quad\quad\quad
\begin{aligned}
R_{319} &=	( xy,\ 1,\ yz) \\
R_{320} &=	( xy+x+y,\ 0,\ yz+y+z) 
\end{aligned}
\end{equation}
\item \begin{equation}
\begin{CD}
R_{321} @>\sigma_1>> R_{322}\\
 @V\sigma_2VV 	 @VV\sigma_2V \\
R_{323} @>>\sigma_1> R_{324}
\end{CD}
\quad\quad\quad
\begin{aligned}
R_{321} &=	( xy,\ x,\ xy+y) \\
R_{322} &=	( yz+y,\ z,\ yz) \\
R_{323} &=	( xy+x+y,\ x,\ xy+x+1) \\
R_{324} &=	( yz+z+1,\ z,\ yz+y+z) 
\end{aligned}
\end{equation}
\item \begin{equation}
\begin{CD}
R_{325}\\
 @V\sigma_2VV \\
R_{326}
\end{CD}
\quad\quad\quad
\begin{aligned}
R_{325} &=	( xy,\ xyz+xy+yz,\ yz) \\
R_{326} &=	( xy+x+y,\ xyz+xz+y,\ yz+y+z) 
\end{aligned}
\end{equation}

\subsection{Solutions of image cardinality 5}
\item \begin{equation}
\begin{CD}
R_{327}\\
 @V\sigma_2VV \\
R_{328}
\end{CD}
\quad\quad\quad
\begin{aligned}
R_{327} &=	( x,\ xyz,\ z) \\
R_{328} &=	( x,\ xyz+xy+xz+yz+x+y+z,\ z) 
\end{aligned}
\end{equation}
\item \begin{equation}
\begin{CD}
R_{329} @>\sigma_1>> R_{330}\\
 @V\sigma_2VV 	 @VV\sigma_2V \\
R_{331} @>>\sigma_1> R_{332}
\end{CD}
\quad\quad\quad
\begin{aligned}
R_{329} &=	( y,\ xz,\ yz+y+z) \\
R_{330} &=	( xy+x+y,\ xz,\ y) \\
R_{331} &=	( y,\ xz+x+z,\ yz) \\
R_{332} &=	( xy,\ xz+x+z,\ y) 
\end{aligned}
\end{equation}
\item \begin{equation}
\begin{CD}
R_{333} @>\sigma_1>> R_{334}\\
 @V\sigma_2VV 	 @VV\sigma_2V \\
R_{335} @>>\sigma_1> R_{336}
\end{CD}
\quad\quad\quad
\begin{aligned}
R_{333} &=	( y,\ xz,\ xyz+xz+yz+y+z) \\
R_{334} &=	( xyz+xy+xz+x+y,\ xz,\ y) \\
R_{335} &=	( y,\ xz+x+z,\ xyz+xy+y) \\
R_{336} &=	( xyz+yz+y,\ xz+x+z,\ y) 
\end{aligned}
\end{equation}
\item \begin{equation}
\begin{CD}
R_{337}\\
 @V\sigma_2VV \\
R_{338}
\end{CD}
\quad\quad\quad
\begin{aligned}
R_{337} &=	( xy,\ y,\ yz) \\
R_{338} &=	( xy+x+y,\ y,\ yz+y+z) 
\end{aligned}
\end{equation}
\item \begin{equation}
\begin{CD}
R_{339} @>\sigma_1>> R_{340}\\
 @V\sigma_2VV 	 @VV\sigma_2V \\
R_{341} @>>\sigma_1> R_{342}
\end{CD}
\quad\quad\quad
\begin{aligned}
R_{339} &=	( xy,\ xz+x+z,\ xy+y) \\
R_{340} &=	( yz+y,\ xz+x+z,\ yz) \\
R_{341} &=	( xy+x+y,\ xz,\ xy+x+1) \\
R_{342} &=	( yz+z+1,\ xz,\ yz+y+z) 
\end{aligned}
\end{equation}

\subsection{Solutions of image cardinality 6}
\item \begin{equation}
\begin{CD}
R_{343} @>\sigma_1>> R_{344}\\
 @V\sigma_2VV 	 @VV\sigma_2V \\
R_{345} @>>\sigma_1> R_{346}
\end{CD}
\quad\quad\quad
\begin{aligned}
R_{343} &=	( x,\ y,\ yz) \\
R_{344} &=	( xy,\ y,\ z) \\
R_{345} &=	( x,\ y,\ yz+y+z) \\
R_{346} &=	( xy+x+y,\ y,\ z) 
\end{aligned}
\end{equation}
\item \begin{equation}
\begin{CD}
R_{347} @>\sigma_1>> R_{348}\\
 @V\sigma_2VV 	 @VV\sigma_2V \\
R_{349} @>>\sigma_1> R_{350}
\end{CD}
\quad\quad\quad
\begin{aligned}
R_{347} &=	( x,\ y,\ yz+y+1) \\
R_{348} &=	( xy+y+1,\ y,\ z) \\
R_{349} &=	( x,\ y,\ yz+z) \\
R_{350} &=	( xy+x,\ y,\ z) 
\end{aligned}
\end{equation}
\item \begin{equation}
\begin{CD}
R_{351} @>\sigma_1>> R_{352}\\
 @V\sigma_2VV 	 @VV\sigma_2V \\
R_{353} @>>\sigma_1> R_{354}
\end{CD}
\quad\quad\quad
\begin{aligned}
R_{351} &=	( x,\ z,\ yz) \\
R_{352} &=	( xy,\ x,\ z) \\
R_{353} &=	( x,\ z,\ yz+y+z) \\
R_{354} &=	( xy+x+y,\ x,\ z) 
\end{aligned}
\end{equation}
\item \begin{equation}
\begin{CD}
R_{355} @>\sigma_1>> R_{356}\\
 @V\sigma_2VV 	 @VV\sigma_2V \\
R_{357} @>>\sigma_1> R_{358}
\end{CD}
\quad\quad\quad
\begin{aligned}
R_{355} &=	( x,\ xy,\ z) \\
R_{356} &=	( x,\ yz,\ z) \\
R_{357} &=	( x,\ xy+x+y,\ z) \\
R_{358} &=	( x,\ yz+y+z,\ z) 
\end{aligned}
\end{equation}
\item \begin{equation}
\begin{CD}
R_{359} @>\sigma_1>> R_{360}\\
 @V\sigma_2VV 	 @VV\sigma_2V \\
R_{361} @>>\sigma_1> R_{362}
\end{CD}
\quad\quad\quad
\begin{aligned}
R_{359} &=	( x,\ xy+x+1,\ z) \\
R_{360} &=	( x,\ yz+z+1,\ z) \\
R_{361} &=	( x,\ xy+y,\ z) \\
R_{362} &=	( x,\ yz+y,\ z) 
\end{aligned}
\end{equation}
\item \begin{equation}
\begin{CD}
R_{363} @>\sigma_1>> R_{364}\\
 @V\sigma_2VV 	 @VV\sigma_2V \\
R_{365} @>>\sigma_1> R_{366}
\end{CD}
\quad\quad\quad
\begin{aligned}
R_{363} &=	( x,\ xy+xz+y,\ z) \\
R_{364} &=	( x,\ xz+yz+y,\ z) \\
R_{365} &=	( x,\ xy+xz+z,\ z) \\
R_{366} &=	( x,\ xz+yz+x,\ z) 
\end{aligned}
\end{equation}
\item \begin{equation}
\begin{CD}
R_{367} @>\sigma_1>> R_{368}\\
 @V\sigma_2VV 	 @VV\sigma_2V \\
R_{369} @>>\sigma_1> R_{370}
\end{CD}
\quad\quad\quad
\begin{aligned}
R_{367} &=	( x,\ yz,\ y) \\
R_{368} &=	( y,\ xy,\ z) \\
R_{369} &=	( x,\ yz+y+z,\ y) \\
R_{370} &=	( y,\ xy+x+y,\ z) 
\end{aligned}
\end{equation}
\item \begin{equation}
\begin{CD}
R_{371} @>\sigma_1>> R_{372}\\
 @V\sigma_2VV 	 @VV\sigma_2V \\
R_{373} @>>\sigma_1> R_{374}
\end{CD}
\quad\quad\quad
\begin{aligned}
R_{371} &=	( x,\ yz,\ yz+y+z) \\
R_{372} &=	( xy+x+y,\ xy,\ z) \\
R_{373} &=	( x,\ yz+y+z,\ yz) \\
R_{374} &=	( xy,\ xy+x+y,\ z) 
\end{aligned}
\end{equation}
\[{\bf R_{371}\colon}\quad \im\delta^2=0
\quad \ker\delta^3=
	\begin{pmatrix}
1 & 0 & 1 & 1 & 0 & -1 & 0 & 0 \\
0 & 1 & 0 & 0 & 0 & 1 & 0 & 0 \\
0 & 0 & 0 & 0 & 1 & 1 & 1 & 1 
	\end{pmatrix}^{\rm T}
\]
\[{\bf R_{372}\colon}\quad \im\delta^2=0
\quad \ker\delta^3=
	\begin{pmatrix}
1 & 0 & 1 & 0 & 0 & -1 & 1 & 0 \\
0 & 1 & 0 & 1 & 0 & 1 & 0 & 1 \\
0 & 0 & 0 & 0 & 1 & 1 & 0 & 0 
	\end{pmatrix}^{\rm T}
\]
\[{\bf R_{373}\colon}\quad \im\delta^2=0
\quad \ker\delta^3=
	\begin{pmatrix}
1 & 1 & 0 & 1 & 0 & 0 & -1 & 0 \\
0 & 0 & 1 & 0 & 0 & 0 & 1 & 0 \\
0 & 0 & 0 & 0 & 1 & 1 & 1 & 1 
	\end{pmatrix}^{\rm T}
\]
\[{\bf R_{374}\colon}\quad \im\delta^2=0
\quad \ker\delta^3=
	\begin{pmatrix}
1 & 0 & 0 & -1 & 1 & 0 & 1 & 0 \\
0 & 1 & 0 & 1 & 0 & 1 & 0 & 1 \\
0 & 0 & 1 & 1 & 0 & 0 & 0 & 0 
	\end{pmatrix}^{\rm T}
\]
\item \label{113} \begin{equation}
\begin{CD}
R_{375} @>\sigma_1>> R_{376}\\
 @V\sigma_2VV 	 @VV\sigma_2V \\
R_{377} @>>\sigma_1> R_{378}
\end{CD}
\quad\quad
\begin{aligned}
R_{375} &=	( xy,\ xz+x+z,\ xyz+xy+y) \\
R_{376} &=	( xyz+yz+y,\ xz+x+z,\ yz) \\
R_{377} &=	( xy+x+y,\ xz,\ xyz+xz+yz+y+z) \\
R_{378} &=	( xyz+xy+xz+x+y,\ xz,\ yz+y+z) 
\end{aligned}
\end{equation}

\subsection{Solutions of image cardinality 7}
\item \begin{equation}
\begin{CD}
R_{379}\\
 @V\sigma_2VV \\
R_{380}
\end{CD}
\quad\quad\quad
\begin{aligned}
R_{379} &=	( x,\ xyz+xz+y,\ z) \\
R_{380} &=	( x,\ xyz+xy+yz,\ z) 
\end{aligned}
\end{equation}

\subsection{Solutions of image cardinality 8}
\item \begin{equation}
\varnothing
\quad\quad\quad
\begin{aligned}
R_{381} &=	( x,\ y,\ z) 
\end{aligned}
\end{equation}
\[{\bf R_{381}\colon}\quad \im\delta^2=0
\quad \ker\delta^3=
	\begin{pmatrix}
1 & 0 & 0 & 0 & 0 & 0 & 0 & 0 \\
0 & 1 & 0 & 0 & 0 & 0 & 0 & 0 \\
0 & 0 & 1 & 0 & 0 & 0 & 0 & 0 \\
0 & 0 & 0 & 1 & 0 & 0 & 0 & 0 \\
0 & 0 & 0 & 0 & 1 & 0 & 0 & 0 \\
0 & 0 & 0 & 0 & 0 & 1 & 0 & 0 \\
0 & 0 & 0 & 0 & 0 & 0 & 1 & 0 \\
0 & 0 & 0 & 0 & 0 & 0 & 0 & 1 
	\end{pmatrix}^{\rm T}
\]
\item \begin{equation}
\begin{CD}
R_{382} @>\sigma_1>> R_{383}
\end{CD}
\quad\quad\quad
\begin{aligned}
R_{382} &=	( x,\ y,\ z+1) \\
R_{383} &=	( x+1,\ y,\ z) 
\end{aligned}
\end{equation}
\begin{multline}
{\bf R_{382}\colon}\quad \im\delta^2=
	\begin{pmatrix}
1 & -1 & 1 & -1 & 1 & -1 & 1 & -1 
	\end{pmatrix}^{\rm T}
\\ \ker\delta^3=
	\begin{pmatrix}
1 & 0 & 0 & -1 & 0 & -1 & 1 & 0 \\
0 & 1 & 0 & 1 & 0 & 1 & 0 & 1 \\
0 & 0 & 1 & 1 & 1 & 1 & 0 & 0 
	\end{pmatrix}^{\rm T}
\end{multline}
\begin{multline}
{\bf R_{383}\colon}\quad \im\delta^2=
	\begin{pmatrix}
1 & 1 & 1 & 1 & -1 & -1 & -1 & -1 
	\end{pmatrix}^{\rm T}
\\ \ker\delta^3=
	\begin{pmatrix}
1 & 0 & 0 & 1 & 0 & -1 & -1 & 0 \\
0 & 1 & 1 & 0 & 0 & 1 & 1 & 0 \\
0 & 0 & 0 & 0 & 1 & 1 & 1 & 1 
	\end{pmatrix}^{\rm T}
\end{multline}
\item \begin{equation}
\varnothing
\quad\quad\quad
\begin{aligned}
R_{384} &=	( x,\ y+1,\ z) 
\end{aligned}
\end{equation}
\begin{multline}
{\bf R_{384}\colon}\quad \im\delta^2=
	\begin{pmatrix}
1 & 1 & -1 & -1 & 1 & 1 & -1 & -1 
	\end{pmatrix}^{\rm T}
\\ \ker\delta^3=
	\begin{pmatrix}
1 & 0 & 0 & 1 & 0 & 1 & 1 & 0 \\
0 & 1 & 0 & -1 & 1 & 0 & -1 & 0 \\
0 & 0 & 1 & 1 & 0 & 0 & 1 & 1 
	\end{pmatrix}^{\rm T}
\end{multline}
\item \begin{equation}
\begin{CD}
R_{385} @>\sigma_1>> R_{386}
\end{CD}
\quad\quad\quad
\begin{aligned}
R_{385} &=	( x,\ y+1,\ z+1) \\
R_{386} &=	( x+1,\ y+1,\ z) 
\end{aligned}
\end{equation}
\begin{multline}
{\bf R_{385}\colon}\quad \im\delta^2=
	\begin{pmatrix}
2 & 0 & 0 & -2 & 2 & 0 & 0 & -2 
	\end{pmatrix}^{\rm T}
\\ \ker\delta^3=
	\begin{pmatrix}
1 & 0 & 0 & 0 & 1 & 0 & 0 & 0 \\
0 & 1 & 1 & 0 & 0 & 1 & 1 & 0 \\
0 & 0 & 0 & 1 & 0 & 0 & 0 & 1 
	\end{pmatrix}^{\rm T}
\end{multline}
\begin{multline}
{\bf R_{386}\colon}\quad \im\delta^2=
	\begin{pmatrix}
2 & 2 & 0 & 0 & 0 & 0 & -2 & -2 
	\end{pmatrix}^{\rm T}
\\ \ker\delta^3=
	\begin{pmatrix}
1 & 1 & 0 & 0 & 0 & 0 & 0 & 0 \\
0 & 0 & 1 & 1 & 1 & 1 & 0 & 0 \\
0 & 0 & 0 & 0 & 0 & 0 & 1 & 1 
	\end{pmatrix}^{\rm T}
\end{multline}
\item \begin{equation}
\begin{CD}
R_{387} @>\sigma_1>> R_{388}
\end{CD}
\quad\quad\quad
\begin{aligned}
R_{387} &=	( x,\ z,\ y) \\
R_{388} &=	( y,\ x,\ z) 
\end{aligned}
\end{equation}
\[{\bf R_{387}\colon}\quad \im\delta^2=0
\quad \ker\delta^3=
	\begin{pmatrix}
1 & 1 & 1 & 1 & 0 & 0 & 0 & 0 \\
0 & 0 & 0 & 0 & 1 & 1 & 1 & 1 
	\end{pmatrix}^{\rm T}
\]
\[{\bf R_{388}\colon}\quad \im\delta^2=0
\quad \ker\delta^3=
	\begin{pmatrix}
1 & 0 & 1 & 0 & 1 & 0 & 1 & 0 \\
0 & 1 & 0 & 1 & 0 & 1 & 0 & 1 
	\end{pmatrix}^{\rm T}
\]
\item \begin{equation}
\begin{CD}
R_{389} @>\sigma_1>> R_{390}\\
 @V\sigma_2VV 	 @VV\sigma_2V \\
R_{391} @>>\sigma_1> R_{392}
\end{CD}
\quad\quad\quad
\begin{aligned}
R_{389} &=	( x,\ x+z,\ x+y) \\
R_{390} &=	( y+z,\ x+z,\ z) \\
R_{391} &=	( x,\ x+z+1,\ x+y+1) \\
R_{392} &=	( y+z+1,\ x+z+1,\ z) 
\end{aligned}
\end{equation}
\begin{multline}
{\bf R_{389}\colon}\quad \im\delta^2=
	\begin{pmatrix}
0 & 0 & 0 & 0 & 2 & 0 & 0 & -2 
	\end{pmatrix}^{\rm T}
\\ \ker\delta^3=
	\begin{pmatrix}
1 & 1 & 1 & 1 & 0 & 1 & 1 & 2 \\
0 & 0 & 0 & 0 & 1 & 0 & 0 & -1 
	\end{pmatrix}^{\rm T}
\end{multline}
\begin{multline}
{\bf R_{390}\colon}\quad \im\delta^2=
	\begin{pmatrix}
0 & 2 & 0 & 0 & 0 & 0 & 0 & -2 
	\end{pmatrix}^{\rm T}
\\ \ker\delta^3=
	\begin{pmatrix}
1 & 0 & 1 & 1 & 1 & 1 & 1 & 2 \\
0 & 1 & 0 & 0 & 0 & 0 & 0 & -1 
	\end{pmatrix}^{\rm T}
\end{multline}
\begin{multline}
{\bf R_{391}\colon}\quad \im\delta^2=
	\begin{pmatrix}
2 & 0 & 0 & -2 & 0 & 0 & 0 & 0 
	\end{pmatrix}^{\rm T}
\\ \ker\delta^3=
	\begin{pmatrix}
1 & 0 & 0 & -1 & 0 & 0 & 0 & 0 \\
0 & 1 & 1 & 2 & 1 & 1 & 1 & 1 
	\end{pmatrix}^{\rm T}
\end{multline}
\begin{multline}
{\bf R_{392}\colon}\quad \im\delta^2=
	\begin{pmatrix}
2 & 0 & 0 & 0 & 0 & 0 & -2 & 0 
	\end{pmatrix}^{\rm T}
\\ \ker\delta^3=
	\begin{pmatrix}
1 & 0 & 0 & 0 & 0 & 0 & -1 & 0 \\
0 & 1 & 1 & 1 & 1 & 1 & 2 & 1 
	\end{pmatrix}^{\rm T}
\end{multline}
\item \begin{equation}
\begin{CD}
R_{393} @>\sigma_1>> R_{394}\\
 @V\sigma_2VV 	 @VV\sigma_2V \\
R_{395} @>>\sigma_1> R_{396}
\end{CD}
\quad\quad\quad
\begin{aligned}
R_{393} &=	( x,\ x+z,\ x+y+1) \\
R_{394} &=	( y+z+1,\ x+z,\ z) \\
R_{395} &=	( x,\ x+z+1,\ x+y) \\
R_{396} &=	( y+z,\ x+z+1,\ z) 
\end{aligned}
\end{equation}
\begin{multline}
{\bf R_{393}\colon}\quad \im\delta^2=
	\begin{pmatrix}
1 & 1 & -1 & -1 & 1 & -1 & 1 & -1 
	\end{pmatrix}^{\rm T}
\\ \ker\delta^3=
	\begin{pmatrix}
1 & 1 & 0 & 0 & 1 & 0 & 1 & 0 \\
0 & 0 & 1 & 1 & 0 & 1 & 0 & 1 
	\end{pmatrix}^{\rm T}
\end{multline}
\begin{multline}
{\bf R_{394}\colon}\quad \im\delta^2=
	\begin{pmatrix}
1 & 1 & -1 & 1 & 1 & -1 & -1 & -1 
	\end{pmatrix}^{\rm T}
\\ \ker\delta^3=
	\begin{pmatrix}
1 & 1 & 0 & 1 & 1 & 0 & 0 & 0 \\
0 & 0 & 1 & 0 & 0 & 1 & 1 & 1 
	\end{pmatrix}^{\rm T}
\end{multline}
\begin{multline}
{\bf R_{395}\colon}\quad \im\delta^2=
	\begin{pmatrix}
1 & -1 & 1 & -1 & 1 & 1 & -1 & -1 
	\end{pmatrix}^{\rm T}
\\ \ker\delta^3=
	\begin{pmatrix}
1 & 0 & 1 & 0 & 1 & 1 & 0 & 0 \\
0 & 1 & 0 & 1 & 0 & 0 & 1 & 1 
	\end{pmatrix}^{\rm T}
\end{multline}
\begin{multline}
{\bf R_{396}\colon}\quad \im\delta^2=
	\begin{pmatrix}
1 & 1 & 1 & -1 & -1 & 1 & -1 & -1 
	\end{pmatrix}^{\rm T}
\\ \ker\delta^3=
	\begin{pmatrix}
1 & 1 & 1 & 0 & 0 & 1 & 0 & 0 \\
0 & 0 & 0 & 1 & 1 & 0 & 1 & 1 
	\end{pmatrix}^{\rm T}
\end{multline}
\item \begin{equation}
\begin{CD}
R_{397}\\
 @V\sigma_2VV \\
R_{398}
\end{CD}
\quad\quad\quad
\begin{aligned}
R_{397} &=	( x,\ xz+y,\ z) \\
R_{398} &=	( x,\ xz+x+y+z+1,\ z) 
\end{aligned}
\end{equation}
\item \begin{equation}
\varnothing
\quad\quad\quad
\begin{aligned}
R_{399} &=	( x+1,\ y,\ z+1) 
\end{aligned}
\end{equation}
\begin{multline}
{\bf R_{399}\colon}\quad \im\delta^2=
	\begin{pmatrix}
2 & 0 & 2 & 0 & 0 & -2 & 0 & -2 
	\end{pmatrix}^{\rm T}
\\ \ker\delta^3=
	\begin{pmatrix}
1 & 0 & 1 & 0 & 0 & -1 & 0 & -1 \\
0 & 1 & 0 & 1 & 1 & 2 & 1 & 2 
	\end{pmatrix}^{\rm T}
\end{multline}
\item \begin{equation}
\varnothing
\quad\quad\quad
\begin{aligned}
R_{400} &=	( x+1,\ y+1,\ z+1) 
\end{aligned}
\end{equation}
\begin{multline}
{\bf R_{400}\colon}\quad \im\delta^2=
	\begin{pmatrix}
3 & 1 & 1 & -1 & 1 & -1 & -1 & -3 
	\end{pmatrix}^{\rm T}
\\ \ker\delta^3=
	\begin{pmatrix}
1 & 0 & 0 & 0 & 0 & 0 & 0 & 1 \\
0 & 1 & 0 & 0 & 0 & 0 & 1 & 0 \\
0 & 0 & 1 & 0 & 0 & 1 & 0 & 0 \\
0 & 0 & 0 & 1 & 0 & 1 & 1 & 3 \\
0 & 0 & 0 & 0 & 1 & -1 & -1 & -3 
	\end{pmatrix}^{\rm T}
\end{multline}
\item \begin{equation}
\begin{CD}
R_{401} @>\sigma_1>> R_{402}
\end{CD}
\quad\quad\quad
\begin{aligned}
R_{401} &=	( x+1,\ z,\ y) \\
R_{402} &=	( y,\ x,\ z+1) 
\end{aligned}
\end{equation}
\begin{multline}
{\bf R_{401}\colon}\quad \im\delta^2=
	\begin{pmatrix}
1 & 1 & 1 & 1 & -1 & -1 & -1 & -1 
	\end{pmatrix}^{\rm T}
\\ \ker\delta^3=
	\begin{pmatrix}
1 & 1 & 1 & 1 & 0 & 0 & 0 & 0 \\
0 & 0 & 0 & 0 & 1 & 1 & 1 & 1 
	\end{pmatrix}^{\rm T}
\end{multline}
\begin{multline}
{\bf R_{402}\colon}\quad \im\delta^2=
	\begin{pmatrix}
1 & -1 & 1 & -1 & 1 & -1 & 1 & -1 
	\end{pmatrix}^{\rm T}
\\ \ker\delta^3=
	\begin{pmatrix}
1 & 0 & 1 & 0 & 1 & 0 & 1 & 0 \\
0 & 1 & 0 & 1 & 0 & 1 & 0 & 1 
	\end{pmatrix}^{\rm T}
\end{multline}
\item \begin{equation}
\begin{CD}
R_{403} @>\sigma_1>> R_{404}
\end{CD}
\quad\quad\quad
\begin{aligned}
R_{403} &=	( y,\ x,\ x+y+z) \\
R_{404} &=	( x+y+z,\ z,\ y) 
\end{aligned}
\end{equation}
\item \begin{equation}
\begin{CD}
R_{405} @>\sigma_1>> R_{406}
\end{CD}
\quad\quad\quad
\begin{aligned}
R_{405} &=	( y,\ x,\ x+y+z+1) \\
R_{406} &=	( x+y+z+1,\ z,\ y) 
\end{aligned}
\end{equation}
\end{enumerate}

\section{Discussion of results}\label{s:discussion}

The usual approach to building (commuting) transfer matrices in statistical physics
requires that the $R$-operators from which they are built must be invertible.
Nevertheless, paper~\cite{K-tm} shows that non-invertible solutions may still generate
interesting algebra involving non-standard transfer matrices, namely, 
``anti-tank-hedgehog'' and ``kagome'' transfer matrices.

The four solutions~(113) of image cardinality six look especially nontrivial
and intriguing.

\subsection*{Acknowledgements}
I would like to thank I.~Korepanov for proposing this work and fruitful discussions.

The work was partially supported by the RFBR grant {\Russian мол\_а № 14-01-31019.}

\end{document}